\documentclass[12pt,reqno]{amsart}
\usepackage{graphicx,color}
\usepackage[centertags]{amsmath}
\usepackage{enumerate}
\usepackage{natbib}

\newcommand{\To}{\longrightarrow}

\begin{document}

\title[Induction and Falsifiability in Statistics]{Broccoli Reduces The Risk of\\\hspace{.1in}\\Splenetic Fever!\\\hspace{.1in}\\{\small The use of induction and falsifiability in \\\hspace{.1in}\\ statistics and model selection}}

\maketitle

\begin{center}
{\bf William M. Briggs}  \\ \vskip .05in General Internal Medicine,
Weill Cornell Medical College \\ 525 E. 68th, Box 46,
New York, NY 10021 \\ \textit{email:} mattstat@gmail.com \\
\vskip .1in
\today
\end{center}

\newpage

\vskip .1in\noindent \textsc{Summary:}
The title, a headline, and a typical one, from a newspaper's ``Health \& Wellness" section, usually written by a reporter who has just read a medical journal, can only be the result of an inductive argument, which is an argument from known contingent premisses to the unknown.  What are the premisses and what is unknown for this headline and what does it mean to statistics?

The importance---and rationality---of inductive arguments and their relation to the frequently invoked, but widely and poorly  misunderstood, notion of `falsifiability' are explained in the context of statistical model selection.  No probability model can be falsified, and no hope for model buidling should be sought in that concept.
\vskip .1in\noindent \textsc{Key words:} Falsifiability; Fisher; Induction; Model complexity; Model selection; Occam's razor; P-values; Popper (Karl); Skill score.
\newpage

\section{Introduction}
Everybody knows that
\begin{equation}
\label{hotp}
	\mbox{Because all the many flames observed before have been hot}
\end{equation}
that this is a good reason to believe 
\begin{equation}
\label{hotc}
	\mbox{that {\it this} flame will be hot.}
\end{equation}
At least, I have never met anybody, regardless of his philosophy, who would be willing to put his hand into a bonfire.  Yet there are philosophers, and statisticians, who will claim that (\ref{hotp}) is {\it not} a good reason to believe (\ref{hotc}), and not only that, but also that there is {\it no} reason to believe (\ref{hotc}).

The argument from (\ref{hotp}) to (\ref{hotc}) is {\it inductive}, which as argument from contingent\footnote{Not logically necessary.} premisses which are, or could have been, observed, to a contingent conclusion about something that has not been, and may not be able to be, observed.  An inductive argument must also have its conclusion say about the unobserved something like what the premisses says about the observed.  The word `like' is sufficiently ambiguous, but this has never troubled philosophers who know an inductive argument ``when they see one" \citep{Sto1982}.

The argument from (\ref{hotp}) to (\ref{hotc}) is also invalid in the strict logical sense that the premiss does not entail the conclusion\footnote{Thus, validity means only that the conclusion is logically entailed by the premisses; invalid does {\it not} imply unreasonable.}.  This should be obvious from the example, because it is possible that the next flame I come upon will not be hot, even though all the other flames I have ever experienced have been.

Regardless of the common sense of (\ref{hotc}), the early part of the 20th century saw the beginning, growth, and dispersal of the belief that {\it all} inductive arguments are unreasonable.  Karl Popper was the philosophical father (Hume was its grandfather) of inductive skepticism.  Thomas Kuhn, Imre Lakatos, Paul Feyerabend and many others are his legtimate children \citep{ThePsi1987}, children who have been increasingly willing to lose touch with reality (this sad history has been recounted in, among other places, \citet{GroLev1994}).   Popper asked, ``Are we rationally justified in reasoning from repeated instances of which we have experience [like the hot flames] to instances of which we have had no experience [this flame]?"  His answer: ``No" \citep{Pop1959}.   This irrational answer has long been exposed for what it is in the analytical philosophical community \citep{Haa2003,ThePsi1987,Sto1982}, but, curiously, the news of its irrationality hasn't reached many scientists yet.  

But, to make the long story mercifully short, Popper convinced himself, and many others, that, since induction could not and should not be trusted, only deduction should be used in scientific inference.  And since it is difficult to prove things positively, Popper therefore claimed that the mark of a `real' theory is that it can be {\it falsified}; theories that could not be were said to be metaphysical and not scientific.  Now, the term {\it falsified} has a precise, unambiguous, logical meaning\footnote{That something was shown to be {\it certainly} false.}, but there are many odd, and incorrect, interpretations of this word in our community, which I will detail below.  

That the falsifiability criterion was nonsensical in the face of theories that are true, and therefore could not be falsified, never bothered Popper.  What did bother him were certain kinds of statistical theories (probability statements) which did not seem to fit into the falsification scheme because, of course, they could not be falsified, and were therefore metaphysical, even though, of course, these theories are in everyday use.  He called this the `problem of decidability,'\footnote{\citet[p. 66]{Sto1982} quotes Hume on this ``custom of calling a {\it difficulty} what pretends to be a {\it demonstration} and endevouring by that means to ellude its force and evidence."} and left it at that; or, rather, he left it for the statisticians to solve \citep{Sto1982}.

Fisher, though certainly not of the same skeptical bent---he often talked about how scientists used inductive reasoning, though he wasn't always entirely clear by what he meant by inductive \citep{Fis1973,Fis1973b}---agreed in principle with the Popperian ideas and used these beliefs to build his system of statistics: theories could only be `rejected' and never verified (and so on).

My purpose of this paper is not to prove that inductive inferences {\it are} reasonable, because that has already been done by others (summarized in \cite{Sto1982}), and, in any case, it is obvious to those of us not infected by the Popperian strain.  I merely want to show that the reasoning behind most statistical methods, and certainly those of model selection, {\it is} inductive, especially when we (civilians and statisticians) step back from the math and try to make sense of what the data tells us. I will also show that falsifiability is of little or no use. These two findings, the importance of induction and the frequently futile search for falsifiability, have important consequences for the future of our field.

\section{Common Inductive Arguments}
About civilians first.  Here is a typical, and schematic, newspaper headline:
\begin{equation}
\label{bp}
	\mbox{Broccoli reduces risk of splenetic fever}
\end{equation}
which had to, of course, come from somewhere.  It is possible, and unfortunately not impossible, that it came directly from the imagination of the reporter.  But it may have also come from an argument like the following:
	\vspace{.1in}
	\begin{center}
		\begin{tabular}{p{3.5in}cc} 
			Broccoli either reduces risk of splenetic fever or it does not & & \\
			&&\\\cline{1-1}
			&&\\
			Broccoli reduces risk of splenetic fever.& & \addtocounter{equation}{1}(\theequation)\\
		\end{tabular}
	\end{center}
	\vspace{.1in}
The premise is a tautology: it is necessarily true regardless of any state of the world.  Now, it is a well known principle of logic that it is impossible to argue from a tautology or necessary truth to a contingent conclusion.   That is, (4) is invalid (and it is not inductive).  It can be made valid if the premiss were changed to: 
$$
\mbox{Broccoli always prevents splenetic fever.}
$$
Evidently, the headline did not come from an argument of this sort, or from the original premiss of (4).  It is more likely that the reporter was reading a medical journal which itself discussed evidence relevant---to the conclusion---from an experiment or an observation on a certain--fixed--group of people.  So the reporter may have been arguing:
	\vspace{.1in}
	\begin{center}
		\begin{tabular}{p{3.5in}cc}
			Of all the people---in this certain group of people---more people who did not eat Broccoli got splenetic fever than did people who ate Broccoli. & & \\
			&&\\\cline{1-1}
			&&\\
			Broccoli reduces risk of splenetic fever.& & \addtocounter{equation}{1}(\theequation)\\
		\end{tabular}
	\end{center}
	\vspace{.1in}
The stated premiss was, obviously, one of the facts reported in the medical journal.  But there are at least two hidden premisses our reporter used, whether he knew them or not: (i) that splenetic fever is unambiguously diagnosed, and (ii) that the facts in the medical journal are accurate.  In any case, I will assume the obvious in what follows: namely, that these, and other similar hidden premisses, are unimportant or do not conflict with the major premisses or conclusion.   

Now, (5) may be valid or invalid depending on whom broccoli reduces the risk and what `reduces risk' means.  If the `whom' is ``for this certain group of people" and `reduces risk' mean ``less people who ate broccoli get splenetic fever", then (5) is valid, but it is merely a tautology and just restates that, in this certain group of people, fewer who ate broccoli got splenetic fever.  However, it is surely false that all newspaper headlines of this sort are tautologies that just repeat the medical findings in different words and that the results hold {\it only} for the certain group of people experimented upon.

Statisticians rarely go to the trouble of tabulating results, like, say, the $2\times 2$ table of broccoli and splenetic fever, and print them in a medical journal without having more than just the certain group of people experimented upon in mind.  That is, they either want to say something, beyond the raw numbers, about specific characteristics of the certain group of people, or about people who are not part of this certain group.   In other words, the conclusion should more specifically state that broccoli reduces risk of splenetic fever for
	\vspace{.1in}
	\begin{center}
		\begin{tabular}{p{3.5in}cc}
			this group such that $P(\mbox{SF}|\mbox{B}) < P(\mbox{SF}|\mbox{No B})$.& & \addtocounter{equation}{1}(\theequation)\\
			&&\\
			or&&\\
			&&\\
			future groups of people not in this certain group.& & \addtocounter{equation}{1}(\theequation)\\
		\end{tabular}
	\end{center}
	\vspace{.1in}
where the abbreviations mean the obvious.  The conclusion (6) says something about this certain group of people, but it says something about an unobservable characteristic of these people; namely the probability of developing splenetic fever given that the people ate, or did not eat, broccoli.  (7) makes a prediction about the presence or absence of splenetic fever for people not in this certain group.  

Either conclusion, (6) or (7), make the argument invalid, but both also make it inductive.  The headline is certainly implying either (6) or (7) or both (it may be implying (6) for future groups of people).  Now, most statistical results, at least those not reported by conscientious statisticians, but certainly those that are found in common medical journals etc., are like this.  That is, it is not clear whether the author's are saying something about unobservable characteristics of their certain group of patients or making predictions about future groups of people.  The implication, I think, for most newspaper's ``Health \& Wellness" sections, is that you should increase your intake of the vegetable/mineral/nutrient-of-the-week mentioned in the headline to certainly avoid developing the disease or malady mentioned.  Meaning, you are to take the headline as a prediction for you.  It is unimportant whether I am right about this, or whether headlines always imply something about unobservables---instead of making predictions about observables in future groups---because most statisticians say something about unobservables.  But it is important that the distinction is hardly ever noted.  

And it's especially critical to understand that whatever the headline means, it is certainly based on an inductive argument.  This is true even if the medical journal's authors were scrupulous in their use of classical statistical methods, and were thus careful to say that it is impossible, based on those methods, to support any  positive conclusion about broccoli and splenetic fever.  Civilians, like our reporter, just do not understand the idiosyncratic and confusing interpretations of p-values and confidence intervals, and they are almost certainly going to go away from a journal believing that the evidence just gathered actually meant something directly about the hypothesis of broccoli reducing the risk of splenetic fever.  Well, so what?  You can argue (incorrectly, I think) that we cannot be responsible for what civilians do with statistics.  But what about statisticians themselves, who {\it do} understand the complexities of classical analysis that know, say, `long-run' is a euphemism for `infinity,' and so on?   What about their arguments?

\section{Popper and Statisticians}
It may be fun to play a game of {\it Who Said It?}:
	\begin{enumerate}[{(}a{)}]
		\item ``We have no reason to believe any proposition about the unobserved {\it even after} experience!"
		
		\item ``There {\it are} no  such things as good positive reasons to believe any scientific theory."
		
		\item ``The truth of any scientific theory is exactly as improbable, both {\it a priori} and in relation to any possible evidence, as the truth of a self-contradictory proposition" (i.e. It is impossible.) 
		
		\item ``Belief, of course, is never rational: it is rational to {\it suspend} belief."
	\end{enumerate}
The first is from the grandfather of inductive skepticism, David Hume \citeyearpar{Hum2003}.  The others are all from Karl Popper \citeyearpar{Pop1959,Pop1963}.  These quotations are important to absorb, because most of us haven't seen them before, and because of {\it that}, a lot of misperceptions about Popper's philosophy and its derivatives are common in our field.  To first show the extent of Popper's influence, we can play another round of our game, this time with current quotes from statisticians:
	\begin{enumerate}[{(}A{)}]	
		\item ```[I]nduction doesn't fit my understanding of scientific (or social scientific) inference." 
		
		\item ``Bayesian inference is good for {\it deductive} inference within a model." {\small (my italics)} 

		\item ``I falsify models all the time."

		\item ``[T]he probability that the `truth' is expressible in the language of probability theory...is vanishingly small, so we should conclude a priori that all theories are falsified." 
		
		\item ``[P]assing such a test does not in itself render [a] theory `proven' or `true' in any sense---indeed, from a thoroughgoing falsificationist standpoint (perhaps even more thoroughgoing than Popper himself would have accepted), we can dispense with such concepts altogether."

		\item ``A theory that makes purportedly meaningful assertions that cannot be falsified by any other observation is `metaphysical.'  Whatever other valuable properties such a theory may have, it would not, in Popper's view, qualify as a {\it scientific} theory."
	\end{enumerate}
It isn't hard to search for examples like this, and there is no reason to hunt for more because these will ring true enough.  The first four quotes are from Andrew Gelman's Columbia statistics blog \citeyearpar{Gel2005}; Dan Navarro wrote the fourth on that blog \citeyearpar{Nav2005}; the last two are from a review paper on Popperism in statistics by \cite[a paper that also contains the line ``Causality does not exist"]{Daw2004}.  (A), I trust, is true, but it is not a statement of logic.  The other comments are, or contain logical statements, and they are false (the second sentence in (F) is a matter of fact and is true).  Before I show that, let me summarize how Popper came to believe what he did, and how these views leaked into statistics.

Hume (it was he who supplied the flames example which started this paper) first to rigorously study the invalidity of inductive inferences \citep{Hum2003}.  Further, he was the first to show that there was no way to remove an inductive argument's invalidity: he proved that no additional, necessarily true or contingent, premisses could be added to the original premisses that would make a given inductive argument valid.  This conclusion is known as {\it inductive fallibilism}, and is nowhere controversial.    

Hume then made an additional step and claimed to have shown that, not only are inductive arguments fallible, but that they were also always unreasonable.   This additional conclusion was shown to hinge on two main premisses: (i) inductive fallibilism, and (ii) {\it deductivism}, which is that all invalid arguments are unreasonable \citep{Sto1982}.  Thus:

	\vspace{.1in}
	\begin{center}
		\begin{tabular}{p{3.5in}cc}
			Inductive fallibilism:  inferences from the observed to unobserved are invalid. & & \\
			&&\\
			Deductivism: all invalid arguments are unreasonable. & & \\
			&&\\\cline{1-1}
			&&\\
			Inductive skepticism:  all inductive arguments are unreasonable. & & \addtocounter{equation}{1}(\theequation)\\
		\end{tabular}
	\end{center}
	\vspace{.1in}

\noindent This is a valid argument, given that both premisses are true.  Again, nobody disputes inductive fallibilism.  How about deductivism?  The flames argument is inductive, therefore invalid, and by deductivism it is {\it unreasonable} to believe that future flames will be hot.  Hume assumed that deductivism was true, but there is no argument that it is; it is taken by him, and by Popper, to be axiomatic.  However, the thesis is plainly wrong \citep{Sto1982}.

Popper took inductive skepticism as his starting point.  Given that the only inferences that are reasonable are deductive ones, and because it is, as mentioned above, impossible to argue from a necessary truth to a contingent conclusion, and all matters of fact are contingent, it becomes impossible to argue directly for the truth of any real-world theory.  The best that you could do is to argue negatively against it: that is, if some theory said that ``$X$" is true, and you directly {\it observed} ``$\neg X$" (not $X$), then you could conclusively say that the theory was false.  But that was all you could do.  You could never say the theory was true, or how likely it was to be true, or whether it was reasonable to believe a theory that was `not yet' falsified, and so on.  That is, Popper argued something like this:

	\vspace{.2in}
	\begin{center}
		\begin{tabular}{p{3.5in}cc}
			(8) & & \\
			&&\\\cline{1-1}
			&&\\
			Theories must be {\it falsifiable} to be `scientific': ``It is a {\it vice} and not a virtue for a model to be infallible." & & \addtocounter{equation}{1}(\theequation)\\
		\end{tabular}
	\end{center}
	\vspace{.2in}

This reasoning, which was fully ``in the air" in the early part of the 20th century, made sense to Fisher, who tried to build falsifiability into his p-values.  Where that got us as a field, by now everybody knows.  However, as I'll show below, and as everybody already knows, p-values cannot falsify a theory; and theories based on probabilty models cannot be falsified, e.g. \cite{Gil1971} and the refutation by \cite{Spi1974} and others.  It may come as a slight surprise to learn that any attempt at using a p-value actually forces its user into making an inductive argument, which are the very things that so horrified Popper.

\section{Induction and Falsifiability in Statistics}
Here are two well-known staples of logic:
	\begin{center}
		\begin{tabular}{cccc}
			 $p\To q$ & & & $p\To q$\\
			$p$& & & $\neg q$\\
			&&&\\\cline{1-1}\cline{4-4}
			&&&\\
			$q$& & & $\neg p$\\
		\end{tabular}
	\end{center}
	\vspace{.1in}
The first, to give it its Latin name and make it official, is {\it modus ponens}, and is to be read ``If (the predicate) $p$ is true, then (the predicate) $q$ is true (or is entailed).  $p$ is true.  Therefore, $q$ is true."  The second, {\it modus tollens}, is to be read ``If $p$ is true, then $q$ is true.  $q$ is false.  Therefore, $p$ is false."  It is these two classic forms, and especially the second, that were latched onto by Popper.  Modus ponens, incidentally, goes from deductive to inductive by replacing the second premiss to ``$q$" (though it doesn't keep its Latin name).

For example, a statistical model (or theory, or hypothesis) $M$ that is truly falsified would have a (valid) argument something like (I take, without further elaboration a `model' $M$ to be the kind of thing that makes statements like ``$M\To q$," where $M$ is not observable; I say nothing here about where models come from):

	\begin{center}
		\begin{tabular}{ccc}
			$M\To P(X>0)=0$& & \\
			$X>0$& & \\
			&&\\\cline{1-1}
			&&\\
			$\neg M$ &&\addtocounter{equation}{1}(\theequation)\\
		\end{tabular}
	\end{center}

\noindent which is to be read ``Model $M$ entails that the probability of seeing an (observable) $X$ greater than 0, is 0; that is, if $M$ is true, it is {\it impossible} that $X>0$.  We saw an $X>0$.  Therefore $M$ is false."  This is great when it happens, as $M$ is {\it deduced} to be false, but this happens rarely in practice, and never does in probabilty models.  Consider instead this more common argument:

	\begin{center}
		\begin{tabular}{ccc}
			$M\To P(X>0)=\epsilon>0$& & \\
			$X>0$& & \\
			&&\\\cline{1-1}
			&&\\
			$\neg M$ && \addtocounter{equation}{1}(\theequation)\\
		\end{tabular}
	\end{center}

\noindent which is to be read ``Model $M$ entails that the probability of seeing an (observable) $X$ is small, as small as you like, but still not zero; that is, it is merely {\it improbable} but {\it not} impossible to see an $X>0$.  We saw an $X>0$ (even a microscopically small $X$).  Therefore $M$ is false."  This argument is not valid, but it is inductive because, of course, no matter how small $P(X>0)$, an $X>0$ might still happen and, when and if it does, it is {\it not} inconsistent with $M$.  It is no good, if you are no fan of induction, rebutting with something on the order of, ``Yes, an $X>0$ is not {\it strictly} inconsistent with $M$, but the probability of seeing such an $X$ given that $M$ is true is so small, that if we do see $X>0$ then $M$ is {\it practically} falsified."  The term `practically falsified' is meaningless and is in the same epistemic boat as `practically a virgin\footnote{I am, it should go without saying, speaking of untouched forestland.}.'  If you insist on something being `nearly' or `practically falsified', then you are making an inductive judgment about $M$, and there is no disguising that fact.  Further, if you choose some cutoff, some particular $\epsilon$, it can be shown that you are also putting a measure of logical probability on the inductive inference for the falsity of $M$ \citep{Jay2003}.

Here is another example which comes close to the silliness of `practically falsified' but is in fact a valid argument: ``[For a series of fair coin flips with M: $P(X_i=H)=0.5$, T]he {\it theoretical event}
$$
n^{-1}\sum_{i=1}^{n}X_i \to 0.5
$$
has M-probability 1.  Hence, as a model of the physical universe, M could be regarded as falsfied if, {\it on observation}, the corresponding physical property, the limiting relative frequency of H in the sequence of coin-tosses exists and equals 0.5, is found to fail" (\citet{Daw2004}; second italics mine; original had `P' instead of `M').  

This argument {\it is} valid, but it is also impossible to fulfill because of the ``on observation" phrase.   Nobody will ever live to see whether the actual limiting frequency of tosses does exist and does falsify $M$\footnote{This was what Keynes was getting at with his ``In the long run we shall all be dead" comment.}.  Stopping at any finite value of tosses, no matter how large, only buys you `practically falsified', which is to say, does not buy you validity, and leaves you holding another inductive inference.

Now, most modern probability models are put into service to say things about unobservable parameters (call them $\theta$).   Here is one possible argument about $M_0$ and $\theta$, where $M_0$ might be a `null' model or hypothesis of some kind, and $\theta>0$ might index some kind of test (say the hypothesis where the mean parameters for two normal groups has that $\theta=\theta_1-\theta_2$; variances known):

	\begin{center}
		\begin{tabular}{ccc}
			$M_0\To P(\theta>0|X)=0$& & \\
			&&\\
			$P(\theta>0|X)=\epsilon>0$& & \\
			&&\\\cline{1-1}
			&&\\
			$\neg M_0$ && \addtocounter{equation}{1}(\theequation)\\
		\end{tabular}
	\end{center}

\noindent The model is to be read, ``If $M_0$ is true, then after I see the data the probability of $\theta>0$ is 0; that is, if $M_0$ is true it is impossible that $\theta>0$.  The actual probability, after seeing $X$, is $P(\theta>0|X)=\epsilon>0$.  Therefore, $M_0$ is false."  This is a valid deductive argument.  Certainly, arguments like this can be made for many, if not all, probability models.  If this kind of argument is what the writer's from (B) and (C) had in mind, then I was wrong and people really are routinely engaged in valid falsifications.  And it may even be true, or `true', as (D) or (E) have it\footnote{Those extraneous marks around the {\it true} are known as {\it scare quotes}, and are a (conscious or not)  attempt by their writers to have it both way about the word in question.  For example: `true' may mean {\it true} or only {\it believed to be true}, which are as far apart in meaning as is possible. It is not clear whether Fisher would have approved of the use of `true' in this sense; we know that he used this technique at least sometimes \citep[p. 334]{Fis1980}  But it is was Popper himself who was the progenitor and true master of this form \citep{Sto1982}. }, that all models are a priori falsified (a claim that actually begins \citet{DenKin2006}).

However, it is clear that conclusions of this type are not what our writers do have in mind.  For, we can reword the conclusion as ``It is false that I am {\it certain} that $\theta>0$; that is, it is false that I know for a fact, without any uncertainty, that $\theta_1>\theta_2$."  So ``$\neg M_0$" merely means ``I am not certain that $\theta_1>\theta_2$",  and that is {\it all} I have gained from this argument; which is to say, I have gained nothing.  

Statisiticians are not interested in models like $M_0$, because probability models start with the tacit assumption that ``I am not certain that $\theta_1>\theta_2$."  This was why an experiment was run and data was collected in the first place.  The tacit assumption is certainly true for the ubiqutous normal model where, no matter what finite set of data is observed, I will never be certain that ``$\theta_1>\theta_2$" is true or false.  The uncertainty is forever built in right at the beginning, and the only way around it is to design a new probability model where, in fact, it is possible to have it certain that $\theta_1>\theta_2$.  But once that is done, it is hard to see how any data would change that fact.

It is also false that all models are a priori falsified.   Presumably, for all observation statements $q$, there is a {\it true} model $M_T$.  It may be, and is even likely, that we will not accurately identify $M_T$.  This does not mean that $M_T$ is falsified, because, of course, it is true.

The best we can do, perhaps, is to identify a set of `useful' models (where I happily leave `useful' vague), none of which are equivalent to $M_T$ (see the discussion in \cite{BerSmi2000}, chap. 6, on ``$M$-closed" vs. ``$M$-complete" vs. ``$M$-open").  It follows that if we knew that these models were {\it not} equivalent to $M_T$, then we would know that the models in the useful set are falsified; in fact, they are {\it all} falsified.  But if we {\it knew} these models were not equivalent to $M_T$, then we would know $M_T$, and it is, again of course, {\it impossible} to falsify what is true (and we wouldn't even bother with creating the useful set, unless we were interested in creating, say, a computationally-simple approximation to $M_T$).

Again, we usually do not know, with certainty, $M_T$.  So we cannot say, with {\it certainty}, that the models in the alternate set are false.  It may be that some models in the set are more useful than others, and to any degree that you like, and this may be all we can ever learn (more on this below).  But they cannot be, a priori or a posteriori, falsified.

Lastly, it worth pointing out that it is not true that we can never know $M_T$, else we could never, for example, create simulations! (also pointed out in \citet[p. 384]{BerSmi2000}).

The classic argument against (but, thanks to Fisher, never {\it for}) a model is:

	\begin{center}
		\begin{tabular}{ccc}
			$M_0\To 0<\mbox{p-value}<1$& & \\
			p-value is small& & \\
			&&\\\cline{1-1}
			&&\\
			$\neg M_0$ && \addtocounter{equation}{1}(\theequation)\\
		\end{tabular}
	\end{center}

\noindent which is to be read,  ``The (null) model $M_0$ entails that we see a uniformly-distributed p-value.  We see a p-value that is publishable (namely, $<0.05$).  Therefore, $M_0$ is false."  This argument is not valid and it is not inductive either because the first premiss says we can see any p-value whatsoever, and since we do (see any value), it is actually evidence {\it for} $M_0$ and not against it. (In fact, if the conclusion were $M_0$, the argument {\it would} be inductive!)  There is {\it no} p-value we could see that would be the logical negation of ``$0<\mbox{p-value }<1$"; well, other than 1 or 0, which may of course happen in practice\footnote{The simplest example is a test for differences in proportion from two groups, where $n_1=n_2=1$ and where $x_1=1, x_2=0$, or $x_1=0, x_2=1$.}.  And when it does, then regardless whether the p-value is 0 {\it or} 1, {\it either} of those values legitimately falsify $M_0$!    

Importantly, the first premiss of (13) is {\it not} that ``If $M_0$ is true, then we expect a `large' p-value,"  because we clearly do not.  But the argument would be valid, and $M_0$ truly falsified, if the first premiss {\it were} ``$M_0\To $ large p-value," but nowhere in the theory of statistics is this kind of statment asserted, though something like it often is.  Fisher was fond of saying---and this is quoted in nearly every introductory textbook---something like this (using my notation):

	\vspace{.1in}
	\begin{center}
		\begin{tabular}{p{3.5in}cc}
			Belief in $M_0$ as an accurate representation of the population sampled is confronted by a logical disjunction: {\it Either} $M_0$ is false, {\it or} the p-value has attained by chance an exceptionally low value \citep[for example]{Fis1970}.  && \\
			&& \addtocounter{equation}{1}(\theequation)\\
		\end{tabular}
	\end{center}
	\vspace{.1in}

\noindent His `logical disjunction' is evidently not one, as the first part of the sentence makes a statement about the unobservable $M_0$, and the second part makes a statement about the observable p-value.  But it is clear that there are implied missing pieces, and his quote can be fixed easily like this:

	\begin{center}
		\begin{tabular}{p{3.5in}cc}
			{\it Either} $M_0$ is false and we see a small p-value, {\it or} $M_0$ is true and we see a small p-value.  
			&& \addtocounter{equation}{1}(\theequation)\\
		\end{tabular}
	\end{center}

\noindent Or just:

	\begin{center}
		\begin{tabular}{p{3.5in}cc}
			{\it Either} $M_0$ is true or it is false and we see a small p-value.  
			&& \addtocounter{equation}{1}(\theequation)\\
		\end{tabular}
	\end{center}

\noindent Since ``{\it Either} $M_0$ is true or it is false" is a tautology, we are left with

	\begin{center}
		\begin{tabular}{p{3.5in}cc}
			We see a small p-value. && \addtocounter{equation}{1}(\theequation)\\
		\end{tabular}
	\end{center}

\noindent Which is of no help at all.  Further, this statement has the same logical status as the a priori judgement in the conclusion of (4); or rather, the p-value casts no direct light on the truth or falsity of $M_0$.  This result should not be surprising, because remember that Fisher argued that the p-value could not deduce whether $M_0$ was true; but if it cannot deduce whether $M_0$ is true, it cannot, logically, deduce whether it is false; that is, it {\it cannot falsify} $M_0$.  

However, current practice is that a small p-value is taken to be by all civilians, and most of us, to mean ``This is evidence that $M_0$ is false."  But that is an inductive argument like this:

	\vspace{.1in}
	\begin{center}
		\begin{tabular}{p{3.5in}cc}
			For most small p-values I have seen in the past, $M_0$ has been false.& & \\
			&&\\
			I see a small p-value and my null hypothesis is $M_0$& & \\
			&&\\\cline{1-1}
			&&\\
			$\neg M_0$ && \addtocounter{equation}{1}(\theequation)\\
		\end{tabular}
	\end{center}
	\vspace{.1in}

This argument has seen success because p-values {\it have} been of some use, but as we know now, it is only because, in simple situations, they are reasonable approximations to (functions of) probability statements of hypotheses like ``$\theta_1>\theta_2$ \mbox{ given } X",  e.g. \citet{BerSel1987}.  

You may also try to salvage (13) by starting with $M_a$ (or with $\neg M_0$), some alternate hypothesis that is not the null hypothesis.  But then, of course, you cannot say anything about a p-value.

\section{Model Selection}
This section is more speculative because I talk about models and how to choose among them.  I do not intend for the statements here to be taken as proof, though I believe that the first two conclusions about the number of models and the existence of a perfect model for any set of data, are true.

How many models are there for any given set of data?  To answer this, \citet{Qui1951,Qui1953} put forth his {\it underdetermination thesis}, which is roughly: for any given model $M$, there will be an indefinite number of other models which are not $M$, but which are equally well supported by the evidence as $M$ is.  This thesis is far from agreed upon \citep{Lis1999,Eng1973,Haa2003}.   But whether or not it is true, it is a fact that people have often used different, non-equivalent, models to explain or predict the same set of observation statements.  Then there is the statement by Kripkenstein that any sequence of numbers has an infinite number of ways the sequence could have been generated \citep{Kri1982,Mad1986}---a thesis which, if true, means that each different way explains {\it and} predicts the observed sequence perfectly.  Again, whether that statement is true, it is again a fact that at least for some sequences, there exists more than one way to generate them.

The conclusion to draw from this is, what may be obvious anyway \citep{BerSmi2000}, that 

	\begin{center}
		\begin{tabular}{p{3.5in}cc}
			There are an infinite number of probability models that can explain any set of data. & & \addtocounter{equation}{1}(\theequation)\\
		\end{tabular}
	\end{center}
	\vspace{.1in}

Now, evidently, for any set of data $x_1,  x_2,  \dots, x_n$, (of any dimensionality) the model $M_\Omega$ exists and says, with a straight face, that, with probability 1, we would have seen just what we saw, namely $x_1$ first, $x_2$ second, and so on (though, conveniently, $M_\Omega$ never reveals itself until after the data comes in: it just always says, unconvincingly, and after the fact, ``I knew it!"\footnote{It was real-life examples of unfalsifiable models like $M_\Omega$ (he mentioned Freudianism and quack medicines as examples) that so (rightly) irritated Popper.}).  So, an argument for $M_\Omega$  might be:

	\begin{center}
		\begin{tabular}{ccc}
			$M_\Omega \To x_1, x_2, \dots, x_n$ & & \\
			$x_1, x_2, \dots, x_n$& & \\
			&&\\\cline{1-1}
			&&\\
			$M_\Omega$ && \addtocounter{equation}{1}(\theequation)\\
		\end{tabular}
	\end{center}

\noindent This is to be read, ``If $M_\Omega$  is true, we will see the data $x_1, x_2, \dots, x_n$, which we do in fact see.  Therefore,  $M_\Omega$ is true."  This argument is not valid, but it is inductive and is some evidence for the truth of $M_\Omega$ in that sense.  It is also an argument, because of the second premiss, only about the already observed data.  It says nothing directly about future $x$s, though it can, of course, be applied to them.  Our experience with such `over-fitted' models can be best stated in the following argument:

	\vspace{.1in}
	\begin{center}
		\begin{tabular}{p{3.5in}cc}
			Of all the many models in the past, simpler ones usually turned out better than complex ones, & & \\
			& & \\
			By (15) there is an $M_\Omega$, and by (14) there is at least an  $M_2\ne M_{\Omega}$, & & \\
			& & \\
			$\mbox{Complexity}(M_{\Omega}) > \mbox{Complexity}(M_2)$. & & \\
			&&\\\cline{1-1}
			&&\\
			$M_2$& & \addtocounter{equation}{1}(\theequation)\\
		\end{tabular}
	\end{center}
	
\noindent This argument is invalid but inductive and is, of course, one version of Occam's Razor.  It is also sufficiently vague because of the terms `better' and `complexity.'  The first term certainly does not mean ``fits the data well", because nothing would ever fit the observed data better than $M_\Omega$, which of course fits without error (where `fit' may be taken in either the `small variance of the parameters' or in the predictive sense).  It may mean ``predicts future data well" (again, it has to be future data, because $M_\Omega$ predicts the present data perfectly).

So, ignore, for a moment, the subject of `complexity' and consider this argument:

	\begin{center}
		\begin{tabular}{p{3.5in}cc}
			$M_2\To \mbox{Score}(M_{\Omega}) < \mbox{Score}(M_2)$ in future data & & \\
			& & \\
			$\mbox{Score}(M_\Omega) < \mbox{Score}(M_2)$ in future data & & \\
			&&\\\cline{1-1}
			&&\\
			$M_2$& & \addtocounter{equation}{1}(\theequation)\\
		\end{tabular}
	\end{center}

\noindent which is to be read, ``If $M_2$ is true, then the prediction score (or negative measure of loss, or utility, or skill, or whatever, but where higher scores are better) for $M_\Omega$ will be less than that for $M_2$.  The score was lower for $M_\Omega$, therefore, $M_2$ is true."  This argument is again invalid and inductive, because no finite set of data, and the score based on them, would insure with certainty that $\mbox{Score}(M_\Omega)$ is always less than $\mbox{Score}(M_2)$.

Very well, suppose `better' in (21) means at least ``predicts future data well."  `Complexity' usually means something like ``effective number of parameters" or ``dimensionality of $\theta$", which are close enough for us here. All this does is change the first premiss of (21) to

	\vspace{.1in}
	\begin{center}
		\begin{tabular}{p{3.5in}cc}
			Of all the many models in the past, ones with fewer (effective) parameters usually predict future data better (give higher scores) than models with more (effective) parameters. & & \addtocounter{equation}{1}(\theequation)\\
		\end{tabular}
	\end{center}
	\vspace{.1in}
\noindent The conclusion remains the same, and the argument is still inductive.  Similar to this, the popular model selection criteria AIC and BIC are, at least partially, based on inductive arguments \citep{Was2000}.

\section{Concluding Remarks}

My conclusion, then, by no means original, is that, in general, there is no formal solution to the problem of model selection.  By `formal', I mean that a procedure that could be followed, in finite time, that would allow the {\it true} model to be deduced, that is, known with certainty, and would allow incorrect ones to be falsified.  Just as important as the lack of formality is (in the Kripkensteinian sense) that there may be some set of (two or more) models that explain {\it and} predict (any set of observations) perfectly: so that the only way to judge between competing models in this set would be to appeal to other, outside criteria (whatever these may be).  

Models are rarely considered in isolation.  When deciding on the truth or falsity of a given model, we often make reference to what this judgment would mean to our belief in other models.  Haack's \citeyearpar{Haa2003} crossword puzzle metaphor about how all models fit together in painting a picture of reality is relevant.  One model supplies the answer to, say, 1 Down, and this answer must be amicable with at least 1 Across, and so on; the size of the puzzle is large and its boundaries somewhat amorphous.  Just how induction fits in with model selection will be discussed in future work.

The arguments used in the course of probabilty modeling and model selection are inductive (mathematical models are found inductively, too \citep{Pol1968}).  But the careful reader will have noticed that nowhere did I attach a probability measure to any of the conclusions of the inductive arguments given above: {\it inductive arguments are not probability statements}.   Probabilities can certainy be found for these conclusions---$p(\theta|x)$ and $p(x_{\mbox{new}}|x_{\mbox{old}})=\int p(x_{\mbox{new}}|\theta,x_{\mbox{old}})p(\theta|x_{\mbox{old}})d\theta$ are common examples.  Of course, all deductive and non-deductive arguments, including inductive ones, are matters of logic, so any probability statements about their conclusions must be statements of logical probability \citep{Jay2003,Key2004}.  This is an undeveloped area in statistics, but it is of fundamental importance, because it is directly applicable to the nature of probability and to what probability models actually say.    

A recent example is a  fascinating paper by \cite{Wag2004} that gives limits of a probabilized version of modes tollens, which gets at what people mean when they say `practically falsified.'  In that paper (and in my notation), he shows that if $p(q|M)=a$ and $p(\neg q)=b$, then $p(\neg M)\to 1$ as $a,b\to 1$, and also as $a,b\to 0$.  Typically, $a=1, 0\le b<1$, and if so, then $b\le p(\neg M)< 1$.  Wagner also shows that these are the best bounds possible.

Falsifiability has also been shown, as it has been in many other places, to be of little use or interest.

\centerline{\textsc{Acknowledgements}} I thank Russell Zaretzki and Rich Levine for discussions which lead to vast improvements in this paper.
%%%%%%
% BIBLIOGRAPHY
\newpage
\bibliographystyle{plainnat}
\bibliography{logic}

\end{document}